\documentclass{gtart}

\def\ifplaintex{\expandafter\ifx\csname documentclass\endcsname\relax}

\ifplaintex 
\else
\fi

\expandafter\ifx\csname beginpicture\endcsname\relax
\expandafter\ifx\csname documentclass\endcsname\relax
\input pictex \else
\input prepictex \input pictex \input postpictex \fi\fi

\def\gt{{\mathsurround=0pt\it $\cal G\mskip-2mu$eometry \&\ 
$\cal T\!\!$opology}}        

\def\gtp{{\mathsurround=0pt\it $\cal G\mskip-2mu$eometry \&\ 
$\cal T\!\!$opology $\cal P\!$ublications}}  

\def\lognumber#1{\def\thelognumber{#1}}
\def\volumenumber#1{\def\thevolumenumber{#1}}
\def\papernumber#1{\def\thepapernumber{#1}}
\def\volumeyear#1{\def\thevolumeyear{#1}}

\def\pagenumbers#1#2{\def\startpage{#1}\def\finishpage{#2}}
\def\published#1{\def\publishdate{#1}}
\def\proposed#1{\def\theproposer{#1}}
\def\seconded#1{\def\theseconders{#1}}
\def\received#1{\def\receiveddate{#1}}
\def\revised#1{\def\reviseddate{#1}}
\def\accepted#1{\def\accepteddate{#1}}

\long\def\asciiabstract#1{\long\def\theasciiabstract{#1}}
\def\asciikeywords#1{\def\theasciikeywords{#1}}

\let\\\par\let\thelognumber\relax
\let\thevolumenumber\relax\let\thepapernumber\relax
\let\thevolumeyear\relax\let\thesamplenumber\relax\let\startpage\relax
\let\finishpage\relax\let\publishdate\relax\let\receiveddate\relax
\let\reviseddate\relax\let\accepteddate\relax\let\theasciititle\relax
\let\theasciiauthors\relax
\let\theasciiabstract\relax\let\theasciikeywords\relax
\let\theasciiemail\relax\let\theshortauthors\relax\let\theshorttitle\relax

\long\def\maketitlep{   

\count0=\startpage

\gt\hfill      
\beginpicture
\setcoordinatesystem units <0.33truein, 0.33truein> point at 2.2 0.9
\setplotsymbol ({$\cal G$})
\plotsymbolspacing=9truept
\circulararc 315 degrees from 0 1 center at 0 0
\setplotsymbol ({$\cal T$})
\circulararc 315 degrees from 1 -1 center at 1 0
\endpicture
\break
{\small\ifx\thesamplenumber\relax 
Volume \else Sample
\fi\thevolumenumber\ (\thevolumeyear)
\startpage--\finishpage\nl
Published: \publishdate}
\vglue 0.5truein plus 0.4fil minus 0.1truein

{\parskip=0pt\leftskip 0pt plus 1fil\def\\{\par\smallskip}{\ifplaintex\large
\else\Large\fi\bf\thetitle}\par\medskip}   

\vglue 0pt plus 0.1fil 

{\parskip=0pt\leftskip 0pt plus 1fil\def\\{\par}{\sc\theauthors}
\par\medskip}

\vglue 0pt plus 0.1fil 

{\small\parskip=0pt\let\newline\\
{\leftskip 0pt plus 1fil\def\\{\par}{\sl\theaddress}\par}
\expandafter\ifx\theemail\relax    
\relax\else\vglue 5pt plus 0.02fil minus 2pt\def\\{\stdspace{\rm 
and}\stdspace} 
\cl{Email:\stdspace\tt\theemail}\fi
\ifx\theurl\relax                  
\relax\else\vglue 5pt plus 0.02fil minus 2pt\def\\{\stdspace{\rm 
and}\stdspace}
\cl{URL:\stdspace\tt\theurl}\fi\par}

\vglue 7pt plus 0.3fil minus 3pt

{\bf Abstract}
\vglue 5pt plus 0.1fil minus 2pt

\theabstract

\vglue 7pt plus 0.3fil minus 3pt

{\bf AMS Classification numbers}\quad Primary:\quad \theprimaryclass

Secondary:\quad \thesecondaryclass

\vglue 5pt plus 0.3fil minus 2pt

{\bf Keywords}\quad \thekeywords

\vglue 10pt plus 0.5fil minus 5pt

{\small  Proposed: \theproposer\hfill Received: \receiveddate\nl
Seconded: \theseconders\hfill 
\ifx\reviseddate\relax                         
Accepted: \accepteddate                        
\else
Revised: \reviseddate                          
\fi}
\eject
}       

\let\maketitlepage\maketitlep
\let\maketitle\maketitlepage

\font\phead=cmsl9 scaled 950
\font=cmsl9 scaled 1050
\font\pnum=cmbx10 scaled 913
\font\lnum=cmbx10 
\font\pfoot=cmsl9 scaled 950
\font=cmsl9 scaled 1050
\ifplaintex
\headline{\vbox to 0pt{\vskip -4.5mm\line{\small\phead\ifnum
\count0=\startpage ISSN 1364-0380 (on line)
1465-3060 (printed) \hfill {\pnum\folio}\else\ifodd\count0\def\\{ }%
\ifx\theshorttitle\relax\thetitle\else\theshorttitle\fi\hfill{\pnum\folio}
\else\def\\{ and }{\pnum\folio}\hfill\ifx\theshortauthors\relax\theauthors
\else\theshortauthors\fi\fi\fi}\vss}}
\footline{\vbox to 0pt{\vglue 0mm\line{\small\pfoot\ifnum\count0=\startpage
\copyright\ \gtp\hfill\else
\gt, Volume \thevolumenumber\ (\thevolumeyear)\hfill\fi}\vss
}}
\else
\makeatletter
\def\@oddhead{{\small\lhead\ifnum\count0=\startpage ISSN 1364-0380 (on line)
1465-3060 (printed) \hfill {\lnum\number\count0}\else\ifodd\count0
\def\\{ }\ifx\theshorttitle\relax \thetitle \else\theshorttitle\fi\hfill
{\lnum\number\count0}\else\def\\{ and }{\lnum\number\count0}
\hfill\ifx\theshortauthors\relax 
\theauthors\else\theshortauthors\fi\fi\fi}}\def\@evenhead{\@oddhead}
\def\@oddfoot{\small\lfoot\ifnum\count0=\startpage\copyright\ \gtp\hfill\else
\gt, Volume \thevolumenumber\ (\thevolumeyear)\hfill\fi}
\def\@evenfoot{\@oddfoot}
\makeatother
\fi

\newwrite\gtoutfile
\long\gdef\makeheadfile{  
{\def\\{, }\def\s{ }
\immediate\openout\gtoutfile head.xxx
\immediate\write\gtoutfile{To: math@arxiv.org}
\immediate\write\gtoutfile{Subject: put or rep NNNNN:pppp}
\immediate\write\gtoutfile{--text follows this line--}
\immediate\write\gtoutfile{Proxy-for: \ifx\theasciiauthors\relax
\theauthors\else\theasciiauthors\fi\s<\ifx\theasciiemail\relax\theemail\else\theasciiemail\fi>}
\immediate\write\gtoutfile{\noexpand\\}
\immediate\write\gtoutfile{Authors: \ifx\theasciiauthors\relax
\theauthors\else\theasciiauthors\fi}
\immediate\write\gtoutfile{Title: \ifx\theasciititle\relax
\thetitle\else\theasciititle\fi}
\immediate\write\gtoutfile{Subj-class: GT or SG or MG etc}
\immediate\write\gtoutfile{MSC-class: \theprimaryclass\ifx\thesecondaryclass\relax\else, \thesecondaryclass\fi}
\immediate\write\gtoutfile{Journal-ref: Geom. Topol. \thevolumenumber
(\thevolumeyear) \startpage-\finishpage}
\immediate\write\gtoutfile{Comments: Published by Geometry and Topology at}
\immediate\write\gtoutfile{\s\s http://www.maths.warwick.ac.uk/gt/GTVol\thevolumenumber/paper\thepapernumber.abs.html}
\immediate\write\gtoutfile{\noexpand\\}
\immediate\write\gtoutfile{}
\ifx\theasciiabstract\relax
\immediate\write\gtoutfile{\theabstract}\else
\immediate\write\gtoutfile{\theasciiabstract}\fi
\immediate\write\gtoutfile{}
\immediate\write\gtoutfile{\noexpand\\}
\immediate\write\gtoutfile{}
\immediate\closeout\gtoutfile}}  

\def\maketitlepage{\maketitlep\makeheadfile}
\let\maketitle\maketitlepage

\def\ifplaintex{\expandafter\ifx\csname documentclass\endcsname\relax}

\ifplaintex 
\else
\fi

\expandafter\ifx\csname beginpicture\endcsname\relax
\expandafter\ifx\csname documentclass\endcsname\relax
\input pictex \else
\input prepictex \input pictex \input postpictex \fi\fi

\def\gt{{\mathsurround=0pt\it $\cal G\mskip-2mu$eometry \&\ 
$\cal T\!\!$opology}}        

\def\gtp{{\mathsurround=0pt\it $\cal G\mskip-2mu$eometry \&\ 
$\cal T\!\!$opology $\cal P\!$ublications}}  

\def\lognumber#1{\def\thelognumber{#1}}
\def\volumenumber#1{\def\thevolumenumber{#1}}
\def\papernumber#1{\def\thepapernumber{#1}}
\def\volumeyear#1{\def\thevolumeyear{#1}}

\def\pagenumbers#1#2{\def\startpage{#1}\def\finishpage{#2}}
\def\published#1{\def\publishdate{#1}}
\def\proposed#1{\def\theproposer{#1}}
\def\seconded#1{\def\theseconders{#1}}
\def\received#1{\def\receiveddate{#1}}
\def\revised#1{\def\reviseddate{#1}}
\def\accepted#1{\def\accepteddate{#1}}

\long\def\asciiabstract#1{\long\def\theasciiabstract{#1}}
\def\asciikeywords#1{\def\theasciikeywords{#1}}

\let\\\par\let\thelognumber\relax
\let\thevolumenumber\relax\let\thepapernumber\relax
\let\thevolumeyear\relax\let\thesamplenumber\relax\let\startpage\relax
\let\finishpage\relax\let\publishdate\relax\let\receiveddate\relax
\let\reviseddate\relax\let\accepteddate\relax\let\theasciititle\relax
\let\theasciiauthors\relax
\let\theasciiabstract\relax\let\theasciikeywords\relax
\let\theasciiemail\relax\let\theshortauthors\relax\let\theshorttitle\relax

\long\def\maketitlep{   

\count0=\startpage

\gt\hfill      
\beginpicture
\setcoordinatesystem units <0.33truein, 0.33truein> point at 2.2 0.9
\setplotsymbol ({$\cal G$})
\plotsymbolspacing=9truept
\circulararc 315 degrees from 0 1 center at 0 0
\setplotsymbol ({$\cal T$})
\circulararc 315 degrees from 1 -1 center at 1 0
\endpicture
\break
{\small\ifx\thesamplenumber\relax 
Volume \else Sample
\fi\thevolumenumber\ (\thevolumeyear)
\startpage--\finishpage\nl
Published: \publishdate}
\vglue 0.5truein plus 0.4fil minus 0.1truein

{\parskip=0pt\leftskip 0pt plus 1fil\def\\{\par\smallskip}{\ifplaintex\large
\else\Large\fi\bf\thetitle}\par\medskip}   

\vglue 0pt plus 0.1fil 

{\parskip=0pt\leftskip 0pt plus 1fil\def\\{\par}{\sc\theauthors}
\par\medskip}

\vglue 0pt plus 0.1fil 

{\small\parskip=0pt\let\newline\\
{\leftskip 0pt plus 1fil\def\\{\par}{\sl\theaddress}\par}
\expandafter\ifx\theemail\relax    
\relax\else\vglue 5pt plus 0.02fil minus 2pt\def\\{\stdspace{\rm 
and}\stdspace} 
\cl{Email:\stdspace\tt\theemail}\fi
\ifx\theurl\relax                  
\relax\else\vglue 5pt plus 0.02fil minus 2pt\def\\{\stdspace{\rm 
and}\stdspace}
\cl{URL:\stdspace\tt\theurl}\fi\par}

\vglue 7pt plus 0.3fil minus 3pt

{\bf Abstract}
\vglue 5pt plus 0.1fil minus 2pt

\theabstract

\vglue 7pt plus 0.3fil minus 3pt

{\bf AMS Classification numbers}\quad Primary:\quad \theprimaryclass

Secondary:\quad \thesecondaryclass

\vglue 5pt plus 0.3fil minus 2pt

{\bf Keywords}\quad \thekeywords

\vglue 10pt plus 0.5fil minus 5pt

{\small  Proposed: \theproposer\hfill Received: \receiveddate\nl
Seconded: \theseconders\hfill 
\ifx\reviseddate\relax                         
Accepted: \accepteddate                        
\else
Revised: \reviseddate                          
\fi}
\eject
}       

\let\maketitlepage\maketitlep
\let\maketitle\maketitlepage

\font\phead=cmsl9 scaled 950
\font=cmsl9 scaled 1050
\font\pnum=cmbx10 scaled 913
\font\lnum=cmbx10 
\font\pfoot=cmsl9 scaled 950
\font=cmsl9 scaled 1050
\ifplaintex
\headline{\vbox to 0pt{\vskip -4.5mm\line{\small\phead\ifnum
\count0=\startpage ISSN 1364-0380 (on line)
1465-3060 (printed) \hfill {\pnum\folio}\else\ifodd\count0\def\\{ }%
\ifx\theshorttitle\relax\thetitle\else\theshorttitle\fi\hfill{\pnum\folio}
\else\def\\{ and }{\pnum\folio}\hfill\ifx\theshortauthors\relax\theauthors
\else\theshortauthors\fi\fi\fi}\vss}}
\footline{\vbox to 0pt{\vglue 0mm\line{\small\pfoot\ifnum\count0=\startpage
\copyright\ \gtp\hfill\else
\gt, Volume \thevolumenumber\ (\thevolumeyear)\hfill\fi}\vss
}}
\else
\makeatletter
\def\@oddhead{{\small\lhead\ifnum\count0=\startpage ISSN 1364-0380 (on line)
1465-3060 (printed) \hfill {\lnum\number\count0}\else\ifodd\count0
\def\\{ }\ifx\theshorttitle\relax \thetitle \else\theshorttitle\fi\hfill
{\lnum\number\count0}\else\def\\{ and }{\lnum\number\count0}
\hfill\ifx\theshortauthors\relax 
\theauthors\else\theshortauthors\fi\fi\fi}}\def\@evenhead{\@oddhead}
\def\@oddfoot{\small\lfoot\ifnum\count0=\startpage\copyright\ \gtp\hfill\else
\gt, Volume \thevolumenumber\ (\thevolumeyear)\hfill\fi}
\def\@evenfoot{\@oddfoot}
\makeatother
\fi

\newwrite\gtoutfile
\long\gdef\makeheadfile{  
{\def\\{, }\def\s{ }
\immediate\openout\gtoutfile head.xxx
\immediate\write\gtoutfile{To: math@arxiv.org}
\immediate\write\gtoutfile{Subject: put or rep NNNNN:pppp}
\immediate\write\gtoutfile{--text follows this line--}
\immediate\write\gtoutfile{Proxy-for: \ifx\theasciiauthors\relax
\theauthors\else\theasciiauthors\fi\s<\ifx\theasciiemail\relax\theemail\else\theasciiemail\fi>}
\immediate\write\gtoutfile{\noexpand\\}
\immediate\write\gtoutfile{Authors: \ifx\theasciiauthors\relax
\theauthors\else\theasciiauthors\fi}
\immediate\write\gtoutfile{Title: \ifx\theasciititle\relax
\thetitle\else\theasciititle\fi}
\immediate\write\gtoutfile{Subj-class: GT or SG or MG etc}
\immediate\write\gtoutfile{MSC-class: \theprimaryclass\ifx\thesecondaryclass\relax\else, \thesecondaryclass\fi}
\immediate\write\gtoutfile{Journal-ref: Geom. Topol. \thevolumenumber
(\thevolumeyear) \startpage-\finishpage}
\immediate\write\gtoutfile{Comments: Published by Geometry and Topology at}
\immediate\write\gtoutfile{\s\s http://www.maths.warwick.ac.uk/gt/GTVol\thevolumenumber/paper\thepapernumber.abs.html}
\immediate\write\gtoutfile{\noexpand\\}
\immediate\write\gtoutfile{}
\ifx\theasciiabstract\relax
\immediate\write\gtoutfile{\theabstract}\else
\immediate\write\gtoutfile{\theasciiabstract}\fi
\immediate\write\gtoutfile{}
\immediate\write\gtoutfile{\noexpand\\}
\immediate\write\gtoutfile{}
\immediate\closeout\gtoutfile}}  

\def\maketitlepage{\maketitlep\makeheadfile}
\let\maketitle\maketitlepage

\def\ifplaintex{\expandafter\ifx\csname documentclass\endcsname\relax}

\ifplaintex 
\else
\fi

\expandafter\ifx\csname beginpicture\endcsname\relax
\expandafter\ifx\csname documentclass\endcsname\relax
\input pictex \else
\input prepictex \input pictex \input postpictex \fi\fi

\def\gt{{\mathsurround=0pt\it $\cal G\mskip-2mu$eometry \&\ 
$\cal T\!\!$opology}}        

\def\gtp{{\mathsurround=0pt\it $\cal G\mskip-2mu$eometry \&\ 
$\cal T\!\!$opology $\cal P\!$ublications}}  

\def\lognumber#1{\def\thelognumber{#1}}
\def\volumenumber#1{\def\thevolumenumber{#1}}
\def\papernumber#1{\def\thepapernumber{#1}}
\def\volumeyear#1{\def\thevolumeyear{#1}}

\def\pagenumbers#1#2{\def\startpage{#1}\def\finishpage{#2}}
\def\published#1{\def\publishdate{#1}}
\def\proposed#1{\def\theproposer{#1}}
\def\seconded#1{\def\theseconders{#1}}
\def\received#1{\def\receiveddate{#1}}
\def\revised#1{\def\reviseddate{#1}}
\def\accepted#1{\def\accepteddate{#1}}

\long\def\asciiabstract#1{\long\def\theasciiabstract{#1}}
\def\asciikeywords#1{\def\theasciikeywords{#1}}

\let\\\par\let\thelognumber\relax
\let\thevolumenumber\relax\let\thepapernumber\relax
\let\thevolumeyear\relax\let\thesamplenumber\relax\let\startpage\relax
\let\finishpage\relax\let\publishdate\relax\let\receiveddate\relax
\let\reviseddate\relax\let\accepteddate\relax\let\theasciititle\relax
\let\theasciiauthors\relax
\let\theasciiabstract\relax\let\theasciikeywords\relax
\let\theasciiemail\relax\let\theshortauthors\relax\let\theshorttitle\relax

\long\def\maketitlep{   

\count0=\startpage

\gt\hfill      
\beginpicture
\setcoordinatesystem units <0.33truein, 0.33truein> point at 2.2 0.9
\setplotsymbol ({$\cal G$})
\plotsymbolspacing=9truept
\circulararc 315 degrees from 0 1 center at 0 0
\setplotsymbol ({$\cal T$})
\circulararc 315 degrees from 1 -1 center at 1 0
\endpicture
\break
{\small\ifx\thesamplenumber\relax 
Volume \else Sample
\fi\thevolumenumber\ (\thevolumeyear)
\startpage--\finishpage\nl
Published: \publishdate}
\vglue 0.5truein plus 0.4fil minus 0.1truein

{\parskip=0pt\leftskip 0pt plus 1fil\def\\{\par\smallskip}{\ifplaintex\large
\else\Large\fi\bf\thetitle}\par\medskip}   

\vglue 0pt plus 0.1fil 

{\parskip=0pt\leftskip 0pt plus 1fil\def\\{\par}{\sc\theauthors}
\par\medskip}

\vglue 0pt plus 0.1fil 

{\small\parskip=0pt\let\newline\\
{\leftskip 0pt plus 1fil\def\\{\par}{\sl\theaddress}\par}
\expandafter\ifx\theemail\relax    
\relax\else\vglue 5pt plus 0.02fil minus 2pt\def\\{\stdspace{\rm 
and}\stdspace} 
\cl{Email:\stdspace\tt\theemail}\fi
\ifx\theurl\relax                  
\relax\else\vglue 5pt plus 0.02fil minus 2pt\def\\{\stdspace{\rm 
and}\stdspace}
\cl{URL:\stdspace\tt\theurl}\fi\par}

\vglue 7pt plus 0.3fil minus 3pt

{\bf Abstract}
\vglue 5pt plus 0.1fil minus 2pt

\theabstract

\vglue 7pt plus 0.3fil minus 3pt

{\bf AMS Classification numbers}\quad Primary:\quad \theprimaryclass

Secondary:\quad \thesecondaryclass

\vglue 5pt plus 0.3fil minus 2pt

{\bf Keywords}\quad \thekeywords

\vglue 10pt plus 0.5fil minus 5pt

{\small  Proposed: \theproposer\hfill Received: \receiveddate\nl
Seconded: \theseconders\hfill 
\ifx\reviseddate\relax                         
Accepted: \accepteddate                        
\else
Revised: \reviseddate                          
\fi}
\eject
}       

\let\maketitlepage\maketitlep
\let\maketitle\maketitlepage

\font\phead=cmsl9 scaled 950
\font=cmsl9 scaled 1050
\font\pnum=cmbx10 scaled 913
\font\lnum=cmbx10 
\font\pfoot=cmsl9 scaled 950
\font=cmsl9 scaled 1050
\ifplaintex
\headline{\vbox to 0pt{\vskip -4.5mm\line{\small\phead\ifnum
\count0=\startpage ISSN 1364-0380 (on line)
1465-3060 (printed) \hfill {\pnum\folio}\else\ifodd\count0\def\\{ }%
\ifx\theshorttitle\relax\thetitle\else\theshorttitle\fi\hfill{\pnum\folio}
\else\def\\{ and }{\pnum\folio}\hfill\ifx\theshortauthors\relax\theauthors
\else\theshortauthors\fi\fi\fi}\vss}}
\footline{\vbox to 0pt{\vglue 0mm\line{\small\pfoot\ifnum\count0=\startpage
\copyright\ \gtp\hfill\else
\gt, Volume \thevolumenumber\ (\thevolumeyear)\hfill\fi}\vss
}}
\else
\makeatletter
\def\@oddhead{{\small\lhead\ifnum\count0=\startpage ISSN 1364-0380 (on line)
1465-3060 (printed) \hfill {\lnum\number\count0}\else\ifodd\count0
\def\\{ }\ifx\theshorttitle\relax \thetitle \else\theshorttitle\fi\hfill
{\lnum\number\count0}\else\def\\{ and }{\lnum\number\count0}
\hfill\ifx\theshortauthors\relax 
\theauthors\else\theshortauthors\fi\fi\fi}}\def\@evenhead{\@oddhead}
\def\@oddfoot{\small\lfoot\ifnum\count0=\startpage\copyright\ \gtp\hfill\else
\gt, Volume \thevolumenumber\ (\thevolumeyear)\hfill\fi}
\def\@evenfoot{\@oddfoot}
\makeatother
\fi

\newwrite\gtoutfile
\long\gdef\makeheadfile{  
{\def\\{, }\def\s{ }
\immediate\openout\gtoutfile head.xxx
\immediate\write\gtoutfile{To: math@arxiv.org}
\immediate\write\gtoutfile{Subject: put or rep NNNNN:pppp}
\immediate\write\gtoutfile{--text follows this line--}
\immediate\write\gtoutfile{Proxy-for: \ifx\theasciiauthors\relax
\theauthors\else\theasciiauthors\fi\s<\ifx\theasciiemail\relax\theemail\else\theasciiemail\fi>}
\immediate\write\gtoutfile{\noexpand\\}
\immediate\write\gtoutfile{Authors: \ifx\theasciiauthors\relax
\theauthors\else\theasciiauthors\fi}
\immediate\write\gtoutfile{Title: \ifx\theasciititle\relax
\thetitle\else\theasciititle\fi}
\immediate\write\gtoutfile{Subj-class: GT or SG or MG etc}
\immediate\write\gtoutfile{MSC-class: \theprimaryclass\ifx\thesecondaryclass\relax\else, \thesecondaryclass\fi}
\immediate\write\gtoutfile{Journal-ref: Geom. Topol. \thevolumenumber
(\thevolumeyear) \startpage-\finishpage}
\immediate\write\gtoutfile{Comments: Published by Geometry and Topology at}
\immediate\write\gtoutfile{\s\s http://www.maths.warwick.ac.uk/gt/GTVol\thevolumenumber/paper\thepapernumber.abs.html}
\immediate\write\gtoutfile{\noexpand\\}
\immediate\write\gtoutfile{}
\ifx\theasciiabstract\relax
\immediate\write\gtoutfile{\theabstract}\else
\immediate\write\gtoutfile{\theasciiabstract}\fi
\immediate\write\gtoutfile{}
\immediate\write\gtoutfile{\noexpand\\}
\immediate\write\gtoutfile{}
\immediate\closeout\gtoutfile}}  

\def\maketitlepage{\maketitlep\makeheadfile}
\let\maketitle\maketitlepage

\def\ifplaintex{\expandafter\ifx\csname documentclass\endcsname\relax}

\ifplaintex 
\else
\fi

\expandafter\ifx\csname beginpicture\endcsname\relax
\expandafter\ifx\csname documentclass\endcsname\relax
\input pictex \else
\input prepictex \input pictex \input postpictex \fi\fi

\def\gt{{\mathsurround=0pt\it $\cal G\mskip-2mu$eometry \&\ 
$\cal T\!\!$opology}}        

\def\gtp{{\mathsurround=0pt\it $\cal G\mskip-2mu$eometry \&\ 
$\cal T\!\!$opology $\cal P\!$ublications}}  

\def\lognumber#1{\def\thelognumber{#1}}
\def\volumenumber#1{\def\thevolumenumber{#1}}
\def\papernumber#1{\def\thepapernumber{#1}}
\def\volumeyear#1{\def\thevolumeyear{#1}}

\def\pagenumbers#1#2{\def\startpage{#1}\def\finishpage{#2}}
\def\published#1{\def\publishdate{#1}}
\def\proposed#1{\def\theproposer{#1}}
\def\seconded#1{\def\theseconders{#1}}
\def\received#1{\def\receiveddate{#1}}
\def\revised#1{\def\reviseddate{#1}}
\def\accepted#1{\def\accepteddate{#1}}

\long\def\asciiabstract#1{\long\def\theasciiabstract{#1}}
\def\asciikeywords#1{\def\theasciikeywords{#1}}

\let\\\par\let\thelognumber\relax
\let\thevolumenumber\relax\let\thepapernumber\relax
\let\thevolumeyear\relax\let\thesamplenumber\relax\let\startpage\relax
\let\finishpage\relax\let\publishdate\relax\let\receiveddate\relax
\let\reviseddate\relax\let\accepteddate\relax\let\theasciititle\relax
\let\theasciiauthors\relax
\let\theasciiabstract\relax\let\theasciikeywords\relax
\let\theasciiemail\relax\let\theshortauthors\relax\let\theshorttitle\relax

\long\def\maketitlep{   

\count0=\startpage

\gt\hfill      
\beginpicture
\setcoordinatesystem units <0.33truein, 0.33truein> point at 2.2 0.9
\setplotsymbol ({$\cal G$})
\plotsymbolspacing=9truept
\circulararc 315 degrees from 0 1 center at 0 0
\setplotsymbol ({$\cal T$})
\circulararc 315 degrees from 1 -1 center at 1 0
\endpicture
\break
{\small\ifx\thesamplenumber\relax 
Volume \else Sample
\fi\thevolumenumber\ (\thevolumeyear)
\startpage--\finishpage\nl
Published: \publishdate}
\vglue 0.5truein plus 0.4fil minus 0.1truein

{\parskip=0pt\leftskip 0pt plus 1fil\def\\{\par\smallskip}{\ifplaintex\large
\else\Large\fi\bf\thetitle}\par\medskip}   

\vglue 0pt plus 0.1fil 

{\parskip=0pt\leftskip 0pt plus 1fil\def\\{\par}{\sc\theauthors}
\par\medskip}

\vglue 0pt plus 0.1fil 

{\small\parskip=0pt\let\newline\\
{\leftskip 0pt plus 1fil\def\\{\par}{\sl\theaddress}\par}
\expandafter\ifx\theemail\relax    
\relax\else\vglue 5pt plus 0.02fil minus 2pt\def\\{\stdspace{\rm 
and}\stdspace} 
\cl{Email:\stdspace\tt\theemail}\fi
\ifx\theurl\relax                  
\relax\else\vglue 5pt plus 0.02fil minus 2pt\def\\{\stdspace{\rm 
and}\stdspace}
\cl{URL:\stdspace\tt\theurl}\fi\par}

\vglue 7pt plus 0.3fil minus 3pt

{\bf Abstract}
\vglue 5pt plus 0.1fil minus 2pt

\theabstract

\vglue 7pt plus 0.3fil minus 3pt

{\bf AMS Classification numbers}\quad Primary:\quad \theprimaryclass

Secondary:\quad \thesecondaryclass

\vglue 5pt plus 0.3fil minus 2pt

{\bf Keywords}\quad \thekeywords

\vglue 10pt plus 0.5fil minus 5pt

{\small  Proposed: \theproposer\hfill Received: \receiveddate\nl
Seconded: \theseconders\hfill 
\ifx\reviseddate\relax                         
Accepted: \accepteddate                        
\else
Revised: \reviseddate                          
\fi}
\eject
}       

\let\maketitlepage\maketitlep
\let\maketitle\maketitlepage

\font\phead=cmsl9 scaled 950
\font=cmsl9 scaled 1050
\font\pnum=cmbx10 scaled 913
\font\lnum=cmbx10 
\font\pfoot=cmsl9 scaled 950
\font=cmsl9 scaled 1050
\ifplaintex
\headline{\vbox to 0pt{\vskip -4.5mm\line{\small\phead\ifnum
\count0=\startpage ISSN 1364-0380 (on line)
1465-3060 (printed) \hfill {\pnum\folio}\else\ifodd\count0\def\\{ }%
\ifx\theshorttitle\relax\thetitle\else\theshorttitle\fi\hfill{\pnum\folio}
\else\def\\{ and }{\pnum\folio}\hfill\ifx\theshortauthors\relax\theauthors
\else\theshortauthors\fi\fi\fi}\vss}}
\footline{\vbox to 0pt{\vglue 0mm\line{\small\pfoot\ifnum\count0=\startpage
\copyright\ \gtp\hfill\else
\gt, Volume \thevolumenumber\ (\thevolumeyear)\hfill\fi}\vss
}}
\else
\makeatletter
\def\@oddhead{{\small\lhead\ifnum\count0=\startpage ISSN 1364-0380 (on line)
1465-3060 (printed) \hfill {\lnum\number\count0}\else\ifodd\count0
\def\\{ }\ifx\theshorttitle\relax \thetitle \else\theshorttitle\fi\hfill
{\lnum\number\count0}\else\def\\{ and }{\lnum\number\count0}
\hfill\ifx\theshortauthors\relax 
\theauthors\else\theshortauthors\fi\fi\fi}}\def\@evenhead{\@oddhead}
\def\@oddfoot{\small\lfoot\ifnum\count0=\startpage\copyright\ \gtp\hfill\else
\gt, Volume \thevolumenumber\ (\thevolumeyear)\hfill\fi}
\def\@evenfoot{\@oddfoot}
\makeatother
\fi

\newwrite\gtoutfile
\long\gdef\makeheadfile{  
{\def\\{, }\def\s{ }
\immediate\openout\gtoutfile head.xxx
\immediate\write\gtoutfile{To: math@arxiv.org}
\immediate\write\gtoutfile{Subject: put or rep NNNNN:pppp}
\immediate\write\gtoutfile{--text follows this line--}
\immediate\write\gtoutfile{Proxy-for: \ifx\theasciiauthors\relax
\theauthors\else\theasciiauthors\fi\s<\ifx\theasciiemail\relax\theemail\else\theasciiemail\fi>}
\immediate\write\gtoutfile{\noexpand\\}
\immediate\write\gtoutfile{Authors: \ifx\theasciiauthors\relax
\theauthors\else\theasciiauthors\fi}
\immediate\write\gtoutfile{Title: \ifx\theasciititle\relax
\thetitle\else\theasciititle\fi}
\immediate\write\gtoutfile{Subj-class: GT or SG or MG etc}
\immediate\write\gtoutfile{MSC-class: \theprimaryclass\ifx\thesecondaryclass\relax\else, \thesecondaryclass\fi}
\immediate\write\gtoutfile{Journal-ref: Geom. Topol. \thevolumenumber
(\thevolumeyear) \startpage-\finishpage}
\immediate\write\gtoutfile{Comments: Published by Geometry and Topology at}
\immediate\write\gtoutfile{\s\s http://www.maths.warwick.ac.uk/gt/GTVol\thevolumenumber/paper\thepapernumber.abs.html}
\immediate\write\gtoutfile{\noexpand\\}
\immediate\write\gtoutfile{}
\ifx\theasciiabstract\relax
\immediate\write\gtoutfile{\theabstract}\else
\immediate\write\gtoutfile{\theasciiabstract}\fi
\immediate\write\gtoutfile{}
\immediate\write\gtoutfile{\noexpand\\}
\immediate\write\gtoutfile{}
\immediate\closeout\gtoutfile}}  

\def\maketitlepage{\maketitlep\makeheadfile}
\let\maketitle\maketitlepage

\lognumber{154}
\volumenumber{5}\papernumber{8}\volumeyear{2001}
\pagenumbers{267}{285}
\accepted{20 March 2001}
\proposed{Ronald Fintushel}
\received{12 December 2000}
\revised{19 February 2001}
\published{24 March 2001}
\seconded{Robion Kirby, Ronald Stern}

\usepackage{amsmath,amsfonts,amscd}
\usepackage{graphicx,psfrag}
\def\psfraga <#1,#2> #3#4{\psfrag #3 
{\smash{\rlap{\kern #1 \raise #2\hbox{#4}}}}}

\def\N{{\mathbb N}}
\def\Z{{\mathbb Z}}

\def\R{{\mathbb R}}

\newtheorem{thr}{Theorem}[section]

\newtheorem{pr}[thr]{Proposition}

\theoremstyle{remark}

\begin{document} 

\title{Homotopy K3's with several symplectic structures}
\author{Stefano Vidussi}
\address{Department of Mathematics, University of 
California\\Irvine, California 92697,USA}
\email{svidussi@math.uci.edu}

\begin{abstract}
In this note we prove that, for any $n \in \N$, there exist a smooth
$4$--manifold, homotopic to a $K3$ surface, defined by applying the
link surgery method of Fintushel--Stern to a certain $2$--component
graph link, which admits $n$ inequivalent symplectic structures.
\end{abstract}

\asciiabstract{
In this note we prove that, for any integer n, there exist a smooth
4-manifold, homotopic to a K3 surface, defined by applying the link
surgery method of Fintushel-Stern to a certain 2-component graph link,
which admits n inequivalent symplectic structures Keywords: Symplectic
topology of 4-manifolds; Seiberg-Witten theory}

\asciikeywords{Symplectic topology, 4-manifolds, Seiberg-Witten theory}


\maketitlepage

\section{Introduction} 
Let $M$ be a smooth, closed, oriented $4$--manifold. 
A {\it symplectic form} $ \omega$ 
on $M$ is a closed 2--form such that $\omega \wedge \omega > 0$. 
To each symplectic form $\omega$ we can associate a 
homotopy class of almost complex structures on $M$, with canonical class 
$K_{\omega} \in H^{2}(M,\Z)$.

 There exist natural equivalence relations among symplectic forms: 
Two symplectic forms $\omega_{0}, \omega_{1}$ on $M$ are said to be  
{\it homotopy equivalent} if there exist a path of symplectic forms 
$\omega_{t}$ connecting them. 
They are said to be {\it equivalent up to diffeomorphism} if there exist a 
diffeomorphism $f \co  M \rightarrow M$ such that $f^{*} \omega_{0} = \omega_{1}$.
We will say that $\omega_{0}, 
\omega_{1}$ are {\it equivalent} if they are connected by a combination of 
homotopies and diffeomorphisms. This amounts to saying that, if we
denote by Sympl$(M)$ the space of symplectic forms on the smooth manifold $M$, 
two symplectic forms are equivalent if and only if they represent 
elements lying in the same connected component of Sympl$(M)$/Diff$(M)$.

 Homotopies preserve the canonical class, while diffeomorphisms 
act by pull-back on it; these properties can be used to distinguish 
symplectic structures.

 McMullen and Taubes \cite{McMT} provided the first example of a smooth
$4$--manifold (together with some generalizations) admitting two inequivalent symplectic structures. The author has then proved, in
\cite{V}, the ``accidental" existence of a third one. These structures are 
distinguished from the combinatorially inequivalent position of their 
canonical class in the convex hull of basic classes. Since the appearance of
that paper, other examples have been pointed out: in \cite{L} some non-simply 
connected complex surfaces (with the reversed orientation) 
are shown to admit two inequivalent structures.

 Two points seem to be worthy of further investigation: the first is to 
exhibit a manifold with inequivalent
structures, having ``simple topology" (among the simply connected examples 
quoted above, the simplest is homeomorphic to the elliptic surface $E(4)$); 
the second is to establish the existence, for any $n$, of a manifold admitting
$n$ inequivalent structures, possibly in a fixed homeomorphism class.

 The aim of this paper is to combine the two points, proving the following
theorem:
\begin{thr} \label{prin} 
For any $n \in \N$ there exist a smooth $4$--manifold, homotopic 
(hence homeomorphic) to a $K3$ surface, admitting $n$ inequivalent symplectic 
forms.\footnote{Ivan Smith kindly informed us of his recent and independent 
construction of a family $X_{n}$, $n \geq 2$, of pairwise non-homeomorphic 
simply connected $4$--manifolds whose $n$-th member admits $n$ inequivalent 
structures distinguished by the divisibility of their canonical class
(I Smith, {\it On Moduli Space 
of Symplectic Forms}, Math. Res. Lett., to appear).} \end{thr}
We outline now the idea of the proof. Following the approach of \cite{McMT} 
(and, implicitly, of \cite{T2}) we start by constructing a family of 
$2$--component links admitting several fibrations over the circle.
Fibered $3$--manifolds are closely related to symplectic $4$--manifolds. 
The fibrations 
partition in equivalence classes (the definition of equivalent fibrations, as 
we will see, follows closely that of symplectic forms).
We will study these equivalence classes by analyzing the Euler classes of the 
fibrations (which play the same role as the canonical classes for symplectic 
structures). Euler classes of fibrations are vertices
of the polyhedron dual to the unit ball of the Thurston norm for the 
link exterior (see \cite{T2}). Different Euler classes correspond to 
non-homotopic fibrations, and the invariance of the Thurston 
(or Alexander) norm 
allows to control equivalence under diffeomorphism.

 Next, we fix a link with (at least) $n$ inequivalent fibrations. 
Using a suitable generalization of the link surgery construction of 
Fintushel--Stern, we will exhibit, in correspondence to each fibration, a 
homotopy $K3$ naturally admitting a symplectic structure, whose canonical class
is related to the Euler class of that fibration. The smooth structure of these 
symplectic manifolds, instead, will not depend on the fibration but only on 
the link. The underlying smooth manifold can be therefore 
endowed of symplectic structures which have different canonical classes,
hence non-homotopic. There is not a natural analogue of the 
Thurston norm in dimension $4$, but there is a good analogue of the Alexander 
norm, related to Seiberg--Witten invariants. We will use the invariance 
of this suitably defined Seiberg--Witten norm, 
in order to discuss the obstructions to the equivalence under diffeomorphism.

 The paper is divided as follows:
Section \ref{classlink} contains the construction of the family of links 
we will use, and the study of their fibrations. The first part of this 
Section is a specialization or a consequence of the general discussion of 
graph links contained in \cite{EN}.
 In Section \ref{conma} we will 
define, for each link, the (generalized) link surgery $4$--manifold, while in 
Section \ref{defsym} we will define, by symplectic fiber sum, the symplectic 
structures on the link surgery manifolds, and discuss their inequivalence.
Section \ref{n3} contains two examples, and some remarks on the
divisibility of the canonical classes.

 On the manifolds considered in this paper 
(as on the previously known examples) 
symplectic structures are distinguished 
by comparing their canonical class. It is a major open question whether there 
exist inequivalent symplectic structures with the same canonical class. A 
negative answer to this question would imply, by standard results on the 
SW basic classes, that the number of inequivalent symplectic structures on a 
fixed smooth manifold is finite.

\section{A class of graph links admitting inequivalent fibrations} 
\label{classlink}
The aim of this Section is to exhibit a family of $2$--component links which 
admit several fibrations over the circle, and whose Alexander and Thurston 
norms can be easily computed. Thanks to the
thorough analysis of graph links, contained in \cite{EN}, we can 
choose among them a family of $2$--component graph links which will give us 
the suitable examples. For reasons that will appear clear later on,
we look for links with odd linking number.

\begin{figure}[ht!] 
\centerline{\small
\psfrag{K}{$K$}
\psfrag{S}{$S$}
\psfrag{H}{$H$}
\psfraga <-1.5pt,0pt> {1}{\scriptsize$1$}
\psfraga <-1pt,0pt> {3}{\scriptsize$3$}
\psfraga <-1pt,0pt> {0}{\scriptsize$0$}
\includegraphics[width=38mm, height=33mm]{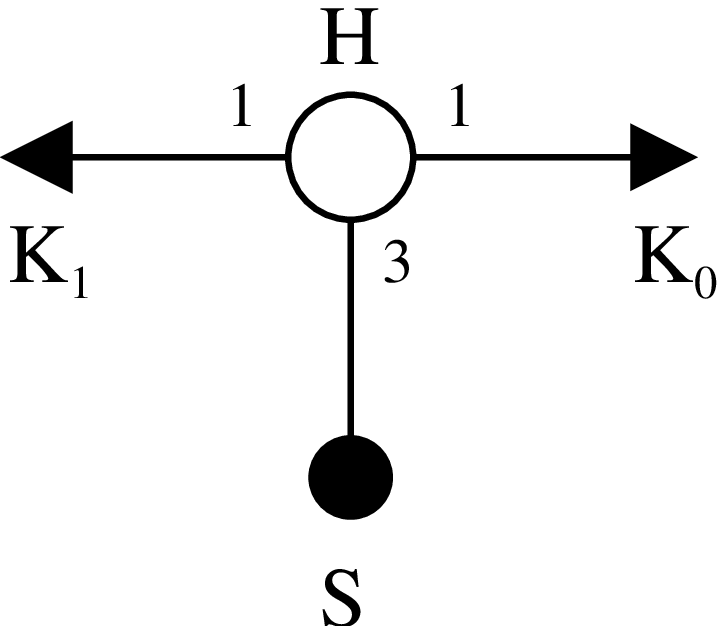}
\hspace*{2cm} 
\includegraphics[width=38mm, height=33mm]{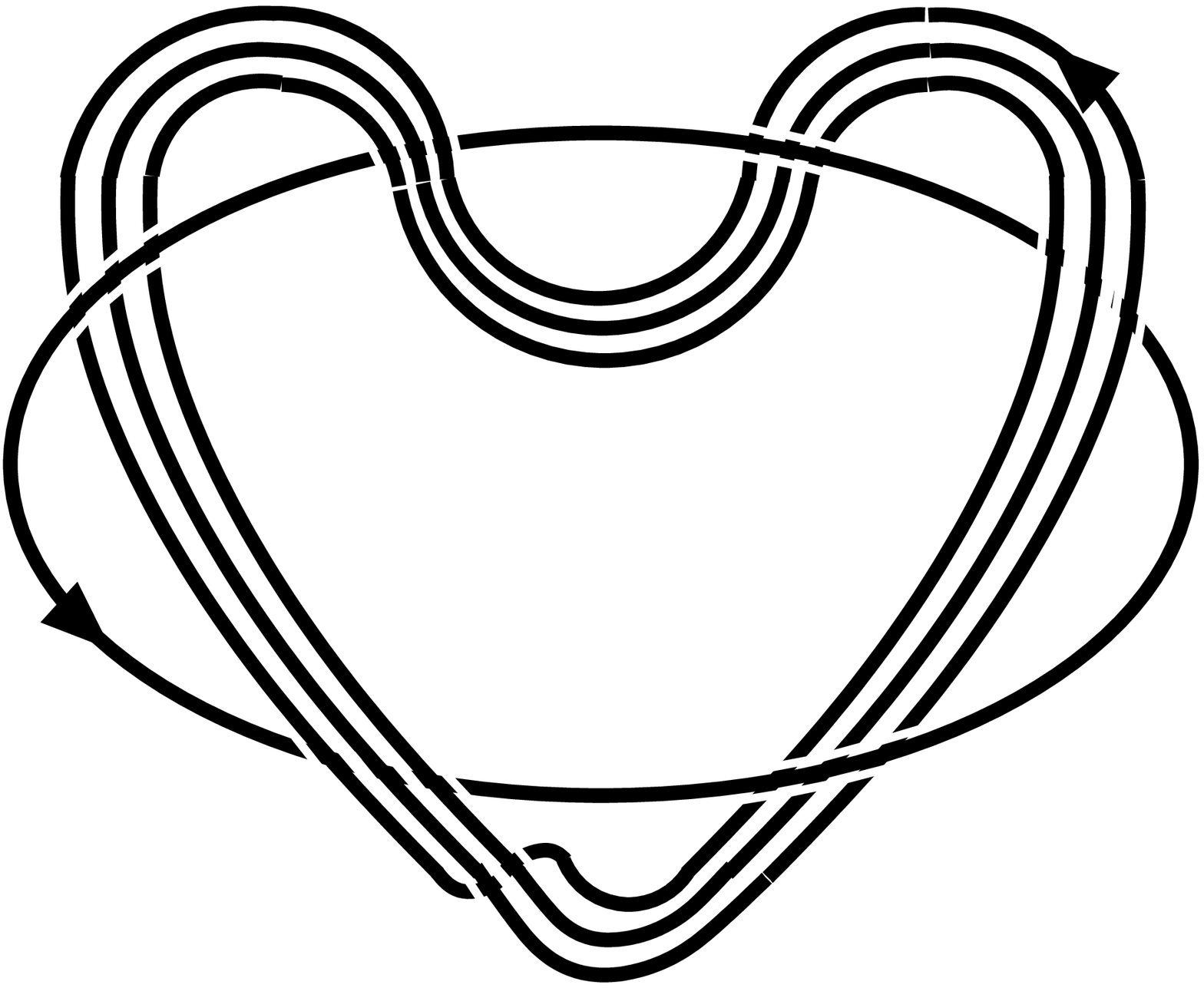}}
\caption{\label{cablin} Diagram representation for $K^{(1)}$ 
(left) and the link $K^{(2)}$ (right)} 
\end{figure}

 Let's introduce the family we will work with. Consider
the Seifert fibration of $S^{3}$ with Seifert invariants $(3,1)$ and take the 
oriented $2$--components link $K^{(1)}$ given
by two regular fibers with their natural orientation. 
This is given by the unknot $K_{0}^{(1)}$ and its $(1,3)$--cable 
$K_{1}^{(1)}$ (order is clearly irrelevant) and is represented, 
in the notation of \cite{EN}, by the diagram on the left in Figure 
\ref{cablin}. 
This diagram requires some explanation. The {\it node} $H$ (a white disk) 
represents the Seifert fibration of the link exterior;
the {\it boundary vertex} $S$ (the black disk) represents the tubular 
neighborhood of the singular fiber of that fibration; the arrowhead vertices 
denote instead the link components $K^{(1)}_{0}$, $K^{(1)}_{1}$. Edge weights 
represent the order of the corresponding fiber: for example, for $K^{(1)}$, the 
arrowheads are regular fibers (order $1$), and the boundary vertex is a 
multiple fiber of order $3$. The node 
(interpreted as regular fiber of the Seifert fibration) and the 
boundary vertex represent, in some sense, {\it virtual} components of 
the link (ie, knots in $S^{3}$ which are not actual components of the link, 
but are naturally associated to its description in terms of the canonical 
decomposition of its exterior).

 Next, take another copy of the same link and splice the two links along the 
$K_{0}$'s (ie, remove a neighborhood of the $K_{0}$'s and glue the link
exteriors identifying meridians with longitudes), obtaining this way a 
$2$--component link (of linking number $9$) in
$S^{3}$. Such an operation correspond to the substitution of the component we 
are splicing with its $(3,1)$--cable (\cite{EN}, Proposition 9.1). 
The image of the link is given on the right of Figure \ref{cablin}. 
Splicing is represented, in diagrams, by joining the edges corresponding to 
the link components we are splicing.

\begin{figure}[ht!] 
\centerline{\small
\psfrag{K}{$K$}
\psfrag{S}{$S$}
\psfrag{H}{$H$}
\psfraga <-1.5pt,0pt> {1}{\scriptsize$1$}
\psfraga <-1pt,0pt> {2}{\scriptsize$2$}
\psfraga <-1pt,0pt> {3}{\scriptsize$3$}
\psfraga <-1pt,0pt> {0}{\scriptsize$0$}
\includegraphics[width=110mm, height=33mm]{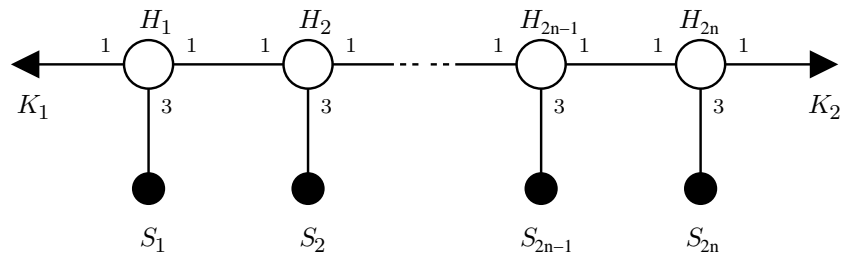}}
\caption{\label{spli}Diagram of the splice decomposition of $K^{(2n)}$} 
\end{figure}

 Iterate $2n$ times
this construction, so to obtain the $2$--component graph link 
$K^{(2n)} \subset S^{3}$. This is a solvable link (link obtained from
unknots by cabling and summing), as the two
components are unknots related by the iterated cabling operations; 
by construction, it has the minimal splice diagram appearing in Figure 
\ref{spli} (minimality refers to the smallest number of edges among
equivalent diagrams, as described in Theorem 8.1 of \cite{EN}).

 In this diagram, each node represents a Seifert fibration appearing in the
JSJ decomposition of the link exterior, the boundary vertices 
represent the tubular neighborhood of a singular fiber of that fibration. 
Together, these represent the virtual components of the link; 
the arrowhead vertices denote instead the {\it real} link components $K_{1}$,
$K_{2}$, in our case regular fibers of the first and last Seifert fibrations. 
The linking number of this link is
$3^{2n}$. This can be computed from the diagram of Figure \ref{spli} 
by multiplying 
the edge weights adjacent to the nodes contained in the path joining $K_{1}$ 
and $K_{2}$ (\cite{EN}, Theorem 10.1). 

 We will denote by
$t_{i}$ the homology class of the meridians to the components $K_{i}$. These 
are the generators of the first homology group of
the link exterior $M^{(2n)} := S^{3} \setminus \nu K^{(2n)}$ (with the symbol 
$\nu ( \cdot )$  we will denote the open 
neighborhood of an embedded submanifold).

 As mentioned above, we are interested in the fibrations of the link 
exterior: we remind the reader that a cohomology class ${\bf m} = (m_{1},m_{2})
\in H^{1}(M,\Z) = [M,S^{1}]$ is said to be fibered (a fibered multilink, in 
the notation of \cite{EN}) if there exist a fibration of $M$ over $S^{1}$ 
whose fibers are minimal genus Seifert surfaces representing ${\bf m}$, ie, 
minimal genus representatives of ${\bf m}$ transverse to the boundary and 
intersecting each boundary component in a cable of the link component, 
whose coefficients and multiplicity are determined by ${\bf m}$ according to 
the formula \begin{equation} \label{boca} 
\sigma_{i} = - (\sum_{j \neq i} m_{j} \cdot lk(K_{j},K_{i}))
\mu(K_{i}) + m_{i} \lambda(K_{i}). \end{equation}
In \cite{EN} the authors have analyzed
in generality the fiberability, the Thurston norm and the Alexander polynomial 
of graph links: if we specialize
their discussion to the case of the link $K^{(2n)}$ we have the 
following result:

\begin{pr}
\label{eisneu} {\rm(a)}\qua Let ${\bf m} := (m_{1},m_{2})$ be a class of 
$H^{1}(M^{(2n)},\Z) = \Z^{2}$; 
then $\bf m$ is fibered if and only if its coordinates
lie in the complement of the intersection of $\Z^{2}$ with the $2n$ lines 
of equation \begin{equation} 3^{i} \cdot m_{1} + 3^{2n - i + 1}
\cdot m_{2} = 0, \ \  \mbox{ $1 \leq i \leq 2n$}; \end{equation} 
as a consequence, the Thurston and the Alexander norms of
$M^{(2n)}$ coincide. 

 {\rm(b)}\qua $\forall \ \ 1 \leq i \leq n$ the Thurston norm 
of the primitive element $P_{i}$ of the $i$-th line, of Cartesian coordinates 
$(3^{2n+1 -2i},-1)$, is given by \begin{equation} 
\label{normp} \|P_{i}\|_{T} = 
3^{4n-2i+1} + 3^{2n-2i+1} + 3^{2n} - 2 \cdot 3^{2n-i} - 2 \cdot 3^{2n-i+1} + 1;
\end{equation} $\forall \ \ n+1 \leq i \leq 2n$ the Thurston norm of the 
primitive element $Q_{i}$ of the $i$-th line, of Cartesian coordinates 
$(1,-3^{2i-1-2n})$, is given by \begin{equation} \label{normq} \|Q_{i}\|_{T} =
3^{2i-1} + 3^{2i-2n-1} + 3^{2n} - 2 \cdot 3^{i} - 2 \cdot 3^{i-1} + 1; 
\end{equation} $\forall \ \ 1 \leq i \leq
n$ the indexed family $a_{i} = \|P_{i}\|_{T}$ is strictly decreasing with 
$i$ and we have the equality \begin{equation}
\|P_{i}\|_{T} = \|Q_{2n+1-i}\|_{T}. \end{equation}

{\rm(c)}\qua The (symmetrized) Alexander polynomial of $K^{(2n)}$ is given by 
\begin{equation} \Delta_{K^{(2n)}}(t_{1},t_{2}) = \prod_{1}^{2n} 
(t_{1}^{3^{i-1}}t_{2}^{3^{2n-i}} + 1 + t_{1}^{-3^{i-1}}t_{2}^{-3^{2n-i}}); 
\end{equation} \end{pr} 

\begin{proof} The proof of this Proposition follows from Theorems 11.1, 11.2
and 12.1 of \cite{EN}.  The input we need for applying these Theorems
is the linking matrix between the link components $K_{j}$ and the
virtual link components corresponding to the regular fiber $H_{i}$ and
the singular fiber $S_{i}$ of each Seifert fibration appearing in the
minimal splice diagram of $K^{(2n)}$, see Figure \ref{spli}.  As we
chose carefully all the Seifert invariants to be equal, the linking
matrix with the $i$-th splice component has the relatively simple form
\begin{equation} \label {linma} \left(
\begin{array}{cc} 
lk(K_{1},H_{i}) & lk(K_{2},H_{i})  \\ \\ lk(K_{1},S_{i}) & lk(K_{2},S_{i}) 
\end{array} \right) = \left( \begin{array}{cc} 3^{i} & 3^{2n-i+1}  
\\ \\ 3^{i-1} & 3^{2n-i} \end{array} \right); \end{equation} the reader 
familiar with the content of \cite{EN} can read this result from the diagram in
Figure \ref{spli} above, as linking numbers coincide 
with the product of the edge 
weights adjacent to the nodes contained in the path from $H_{i}$ (or
$S_{i}$) to $K_{j}$ but not on the path itself.

 To prove part (a), we recall that
Theorem 11.2 of \cite{EN} claims that ${\bf m}$ is fibered if and only if 
its coefficients satisfy 
\begin{equation} \sum_{j = 1}^{2} m_{j} \cdot lk(K_{j},H_{i}) 
\neq 0; \ \ \sum_{j = 1}^{2} m_{j} \cdot lk(K_{j},S_{i}) \neq 0, \ \ 
\forall \ \ 1 \leq i \leq 2n. \end{equation} Substituting with the result
of Equation \ref{linma} we obtain the statement in
(a) above. All the $4n$ (top dimensional) faces of the Thurston 
unit ball are consequently fibered, ie, all the integral point laying in the 
cone over these faces are represented by fibrations.
It is well known (see eg \cite{McM}) that in this case the Thurston and 
Alexander norm coincide. 

 Concerning (b), in Theorem 11.1 of \cite{EN} 
it is proven the Thurston norm of any element ${\bf m} = (m_{1},m_{2})$ 
of $H^{1}(M^{(2n)},\Z)$ is given by 
$$\|(m_{1},m_{2}) \|_{T} =\kern4in$$ \vspace{-6mm}
\begin{equation}\sum_{v=1}^{2n} (\delta_{v} - 2) |  \sum_{j=1}^{2} m_{j} \cdot lk(K_{j},H_{v})|
+ \sum_{w=1}^{2n} (\delta_{w} - 2) |  \sum_{j=1}^{2} m_{j} \cdot 
lk(K_{j},S_{w})| \end{equation} 
where 
$\delta_{( \cdot )}$ is the number of incident edges of the virtual 
component (in our case, $3$ for the nodes, $1$ for the boundary vertices). 
Plugging in this formula the coordinates of $P_{i}, Q_{i}$, and using the 
values of the linking matrix of Equation \ref{linma} we obtain Equations \ref{normp} 
and \ref{normq}. 
The rest of (b) follows from these results. Figure \ref{piqu} gives the 
example of the Thurston unit ball for $n=2$.

 Finally, Theorem 12.2 of \cite{EN} computes the Alexander polynomial of the
link as\newpage
$$\Delta_{K^{(2n)}}(t_{1},t_{2}) =\kern3.8in$$ \vspace{-6mm}
\begin{equation}
\prod_{v=1}^{2n}(t_{1}^{lk(K_{1},H_{v})} t_{2}^{lk(K_{2},H_{v})} -1)^{\delta_{v} - 2} \cdot \prod_{w=1}^{2n}(t_{1}^{lk(K_{1},S_{w})} t_{2}^{lk(K_{2},S_{w})} -1)^{\delta_{w} - 2};
\end{equation} by using Equation \ref{linma} again and taking the 
symmetrized polynomial we obtain (c). \end{proof}

Notice that the Alexander polynomial is invariant with respect to the 
interchange of the two variables. This fact, together with the Blanchfield 
duality 
$\Delta_{K^{(2n)}}(t_{1}^{-1},t_{2}^{-1}) = \Delta_{K^{(2n)}}
(t_{1},t_{2})$ is reflected in the symmetry of the non-fibered lines 
and of the norm of their primitive elements
about the diagonals $m_{1} = m_{2}, m_{1} = - m_{2}$ 
(see eg Figure \ref{piqu}).

 The three dimensional scene is almost complete now. 
We need to remind the reader (see the Introduction of \cite{McMT} for example) 
that to any real cohomology class of $M^{(2n)}$ lying on the cone over each of the $4n$ fibered faces of the Thurston unit ball we can associate a closed, 
nowhere vanishing representative. The classes that are image of integral 
cohomology classes are represented by the fibrations of $M^{(2n)}$ considered 
above.
A fibered face is identified by an element of $H_{1}(M^{(2n)},\Z)$, which is 
the Euler class of the (measured) foliation defined by closed, nowhere 
vanishing forms. This class (up to a factor $-1$) is the vertex 
dual to the fibered face in the dual Thurston polyhedron.

 We can define equivalence classes of closed, nowhere vanishing forms in a 
fashion similar to the one used for symplectic forms:
Two forms are said to be {\it equivalent} if they are
connected by a combination of homotopies (through non-vanishing closed forms) 
and diffeomorphisms of $M^{(2n)}$.

 We can use the Euler class to study the equivalence classes of forms on 
$M^{(2n)}$. Although we will not use directly this equivalence result, 
the method used in 
the proof provides the model for the proof of Theorem \ref{prin}. We have:
\begin{thr} \label{ineq} The moduli space of closed non-vanishing
forms on $M^{(2n)}$ contains at least $n+1$ components. \end{thr}

\begin{proof}
We start by observing that homotopic forms must lie in the cone over the same 
fibered face of the Thurston unit ball.
We will prove the Theorem by showing that the action of 
Diff$(M^{(2n)})$ on $H_{1}(M^{(2n)},\Z)$ partitions the $4n$ Euler classes
in $(n+1)$ different orbits.

 The diffeomorphisms of $M^{(2n)}$ act over
$H^{1}(M^{(2n)},\R)$ through isometries of the Thurston (Alexander) norm, 
preserving the unit ball; they preserve moreover the lattice of integral 
points (and consequently, divisibility).

 We will leave aside the case of $n = 1$, where the statement of the Theorem 
holds true by direct check of the divisibility of the vertices dual to the 
fibered faces ($t_{1}^{8}t_{2}^{8}$ and $t_{1}^{4}t_{4}^{-4}$ are two 
inequivalent vertices, that we can easily determine from
the Alexander polynomial).
\begin{figure}[ht!] 
\centerline{\small
\psfrag{P}{$P$}
\psfrag{Q}{$Q$}
\psfrag {-P}{--$P$}
\psfrag {-Q}{--$Q$}
\psfraga <0.5pt,0pt> {1}{\scriptsize$1$}
\psfraga <0.5pt,0pt> {2}{\scriptsize$2$}
\psfraga <0.5pt,0pt> {3}{\scriptsize$3$}
\psfraga <0.5pt,0pt> {4}{\scriptsize$4$}
\includegraphics[width=60mm, height=55mm]{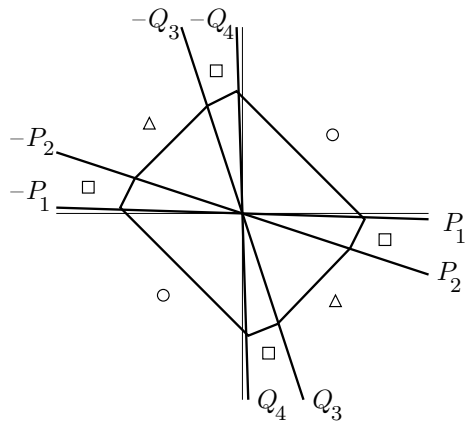}}
\caption{\label{piqu}Thurston unit ball for $K^{(4)}$ --- non-fibered rays
are marked by their primitive element.} 
\end{figure}

 Let's assume then that $n \geq 2$; any primitive element $\pm P_{i},
\pm Q_{i}$, lying in one non-fibered ray must be sent to a 
primitive element lying in a non-fibered ray and having the same norm, ie, 
Diff$(M^{(2n)})$ acts on the set of $4n$ elements 
$\{ \pm P_{j}; \pm Q_{n+j},j=1,...,n  \}$ preserving the norm. 
Proposition \ref{eisneu}, (b) implies that there are 
at least $n$ different orbits for this action: the orbit of $P_{i}$ can contain
only $\pm P_{i}, \pm Q_{2n+1-i}$, as these are the only primitive non-fibered
integral homology classes having the same Thurston norm, and the same
holds for the other $(4n-1)$ elements.
As $n \geq 2$, the knowledge of the action of 
Diff$(M^{(2n)})$ on a single primitive element determines completely the 
action on the others (eg, if $P_{1}$ goes to $-Q_{2n}$, then $P_{2}$ must go 
to $-Q_{2n-1}$ and so on). This implies a constraint on the possible orbit 
of fibered faces, and it is easy to check that there must be at least 
$(n+1)$ orbits of the action of Diff$(M^{(2n)})$ on the set
of fibered faces, the possible identifications corresponding to reflections 
through the diagonals of
$H^{1}(M^{(2n)},\R)$ (the `+1' comes from the fact that the fibered faces 
hit by the diagonals have orbits containing at most $2$ elements, while the 
orbits of the remaining faces contain at most $4$ elements).
Figure \ref{piqu} illustrates the case of $n=2$: the fibered faces marked 
with a different symbol in Figure are inequivalent.
This result implies that the vertices dual to inequivalent fibered faces 
(and consequently the Euler classes) lie in different orbits of 
Diff$(M^{(2n)})$, and completes the proof. \end{proof}
In fact, the existence of inequivalent fibered faces already implies the proof 
of the Theorem, but as in dimension $4$ we will have to work directly with the 
canonical classes rather than the symplectic forms, we want to emphasize
the role of the Euler class. 

 We want to point out that the idea
applied in \cite{McMT} for proving the inequivalence of the Euler classes,
namely the study of the valence of the vertices of the dual Thurston 
polyhedron, cannot be applied here (nor in the proof of the following
Theorem \ref{main}).

 It is conceivable to use direct
informations on the Euler classes to deduce their non-equivalence, for example their divisibility. In some specific examples (with small value of
$n$) it is easy to do so, but for the general case this does not seem to be 
the simpler approach. Note that to determine explicitly the Euler classes
we need only information contained in the
Alexander polynomial, to determine the vertices of the convex hull of the 
Alexander  ``basic classes" (homology
classes with non-zero coefficient in $\Delta_{K^{(2n)}}(t_{1},t_{2})$). 
Or, perhaps more simply, we can observe that the 
formulae contained in (a) 
and (b) of Proposition \ref{eisneu} allow to
determine completely the vertices of the convex hull 
even without recurring to the Alexander polynomial.
We will comment further the question of divisibility in Section \ref{n3}.

 We finish this Section by pointing out that the proof of Theorem 
\ref{ineq} above holds using both the Thurston norm or the
Alexander norm: we will use an analog of the latter 
in the proof of Theorem \ref{prin}.

\section{Construction of the 4--manifolds} \label{conma}
We recall the definition
of a link surgery manifold in this context (for the general case and the details, see \cite{FS}): Consider a
$2$--component oriented link $K \subset S^{3}$ and let $M := S^{3} \setminus \nu K$. Choose a basis
$(\alpha_{i},\beta_{i})$ of simple curves, having intersection $1$, in the homology of the two tori that constitute the boundary of $M$. Next take $2$
copies of the rational elliptic surface
$E(1)$, containing each one an elliptic fiber $F_{i}$, and define
\begin{equation} \label{usdef} E(2,K) = (\coprod E(1)_{i} \setminus \nu F_{i}) \cup_{F_{i} \times S^{1} =
S^{1} \times
\alpha_{i} \times \beta_{i}} S^{1} \times M \end{equation} where
the gluing is made lifting a diffeomorphism between $S^{1} \times
\alpha_{i}$ and $F_{i}$ to an orientation reversing
diffeomorphism of the boundary tori in such a way
that the homology class of $\beta_{i}$ is sent over the homology class
of the normal circle to the $i$-th elliptic fiber. There are usually ambiguities, in the definition of link
surgery manifolds, related to the choice of the $(\alpha_{i},\beta_{i})$, the homology basis for the elliptic
fiber and the lifting of the diffeomorphism, but in this case, as any orientation preserving diffeomorphism of
$\partial (E(1) \setminus \nu F)$ extends to an orientation preserving
diffeomorphisms of $(E(1) \setminus \nu F)$ (see \cite{GS}), the manifold is well defined and moreover its smooth
structure is unaffected by the choice of the homology basis for $\partial M$. The SW
polynomial of $E(2,K)$ has been computed in \cite{FS} and we have \begin{equation} \label{swpol} SW_{E(2,K)} = \pm \Delta_{K}(t_{1}^{2},t_{2}^{2}), 
\end{equation} where $\Delta_{K}$ is the symmetrized Alexander
polynomial of the link, and $t_{i}$ is the image of the class of the 
meridian $\mu(K_{i})$ (that we denoted in the same way) under the
injective map \begin{equation} \label{injmap} H_{1}(M,\Z) \stackrel{S^{1} \times ( \cdot
)}{\rightarrow}  H_{2}(S^{1} \times M,\Z) \rightarrow H_{2}(E(2,K),\Z)
\stackrel{PD}{\rightarrow} H^{2}(E(2,K),\Z). \end{equation} The manifold defined above has the rational homotopy
type of $E(2) = K3$. The Blanchfield duality of the multivariable Alexander polynomial, applied to the two
component link, shows that whenever $K$ has odd linking number, then all basic classes are even, so that $E(2,K)$
is a homotopy $K3$ manifold (in the case at hand, this can also be checked from direct inspection of the $SW$
polynomial). 

 The family of homotopy $K3$'s that will allow us to prove our main result is given, in the smooth
category, by $E(2,K^{(2n)})$. We want to point out that these manifolds are not new, in the sense that they are
contained in the collection of exotic $K3$'s of Fintushel--Stern; their 
definition in \cite{FS} corresponds to the
choice of $(\alpha_{i}, \beta_{i}) = (\mu(K_{i}), - 3^{2n} \mu(K_{i}) + \lambda(K_{i}))$, as $lk(K_{1},K_{2}) =
3^{2n}$, but any choice gives the same smooth manifold.

\section{Definition of the symplectic structures} \label{defsym}
In this Section we will show how we can induce, from inequivalent fibrations 
on $M^{(2n)}$,
inequivalent symplectic structures to $E(2,K^{(2n)})$. We want to warn the reader that, although inspired by
the same bottom line, the idea we will use to construct our examples is slightly different from the one on
\cite{McMT} (in particular, it depends on the use of $E(1)$). 
We will comment on this at the end of the Section.

 For sake of notation, we will
omit the upper index for the link exteriors $M^{(2n)}$ whenever unnecessary.

 Take a fibered cohomology class ${\bf m} \in H^{1}(M,\Z)$.
Recall that such a choice induces a choice (up to isotopy)
of a basis element $\beta_{i}$ for the boundary tori of $M$. 
This is a simple curve determined by the curve $\sigma_{i}$ described in 
Equation \ref{boca} from the relation $\sigma_{i} = d_{i} \beta_{i}$, 
where $d_{i} > 0$ is the 
divisibility of $\sigma_{i}$. 
Choose next any simple curve $\alpha_{i}$ in the boundary tori that 
completes $\beta_{i}$ to a basis, as required in Section \ref{conma}.
The definition of \cite{FS}, as seen in the previous Section, corresponds to 
the class ${\bf m} = (1,1)$ (so in particular we have necessarily $d_{i} = 1$) 
and $\alpha_{i} = \mu(K_{i})$.

 We can associate to ${\bf m}$ a closed three manifold $M({\bf m})$ by
Dehn filling the link exterior $M$ with two solid tori $S^{1} \times D^{2}$, 
using as slope (attaching circle) the two curves $\beta_{i}$. It is a well 
known fact that this prescription is enough to describe, up to diffeomorphism,
the three manifold. This manifold has first Betti number equal to $1$. 
The fibration $\bf m$ of $M$ over $S^{1}$ extends to the closed manifold,
through the fibration of degree $d_{i}$ of $S^{1} \times D^{2}$ over $S^{1}$ on each solid torus. The cores
$C_{i}$ of the solid tori (ie, the so-called dual link in $M({\bf m})$) 
come with a natural framing and clearly
satisfy \begin{equation} \label{dual} M = S^{3} \setminus \nu(K_{1}
\cup K_{2}) = M({\bf m}) \setminus \nu(C_{1} \cup C_{2}). \end{equation} Up to isotopy, these cores are
transverse to the closed fibers of $M({\bf m})$. 

 It is a Theorem of Thurston \cite{T1} that $S^{1} \times
M({\bf m})$ admits a symplectic structure; the embedded, $0$ self-intersection tori $S^{1} \times C_{i}$
are symplectic submanifolds of $S^{1} \times M({\bf m})$ (although, if some $d_{i} \neq 1$, they are not
sections of the natural fibration over $T^{2}$), and come with a trivialization of their normal bundle. We use
these submanifolds to present the manifolds
$E(2,K)$ as fiber sum: we have the following:

 \begin{pr} \label{fibsum} For any choice of fibered ${\bf m} \in
H^{1}(M,\Z)$ we have \begin{equation} \label{ncs} E(2,K) = (\coprod E(1)_{i} ) \#_{F_{i} = S^{1} \times C_{i}}
S^{1} \times M({\bf m}). \end{equation} \end{pr}

\begin{proof} This is a consequence, as in Section \ref{conma}, 
of the property of diffeomorphism extension of
$E(1)$; first, the fiber sum in Equation \ref{ncs} above is well defined, and
independent of the choice of the framing for $S^{1} \times C_{i}$; 
then, because of the diffeomorphism  
of Equation \ref{dual}, for all choice of ${\bf m}$ 
the fiber sum manifold is defined as one of those in Equation \ref{usdef}, where 
we determine the pair of $\beta_{i}$ according to Equation
\ref{boca} and we complete it to a homology basis of $\partial M$ respecting 
the framing of $S^{1} \times
C_{i}$. Up to diffeomorphism, the smooth manifold defined this way
is the same for all ${\bf m}$.
\end{proof} The interest of Proposition
\ref{fibsum} above stems from the fact that, whenever we perform the fiber sum
along symplectic tori, the resulting manifolds admit a symplectic structure.

 Our goal is to discuss the equivalence of the
symplectic structures arising from different 
choices of ${\bf m}$, as we did for the fibrations in Section \ref{classlink}.
In order to do so, we need a norm on homology classes of $E(2,K)$ 
that plays the same role of the Thurston (Alexander) norm on $M$.
We define the Seiberg--Witten norm on $H_{2}(E(2,K),\R) = \R^{b}$
(with $b = 22$) as \begin{equation} \| \varphi \|_{SW} =  \mbox{max} \{ 
k_{i}(\varphi), \ \ \mbox{$k_{i}$ basic class of $E(2,K)$} \}; \end{equation}
a description of the relation between the SW basic classes and the 
SW polynomial can be found, eg, in \cite{FS}.
To study this norm we find convenient to choose a cohomology basis for
$H^{2}(E(2,K),\Z)$ by taking the first two generators to be 
$t_{1},t_{2}$ defined as in Section \ref{conma} and completing them to a basis 
$t_{1},t_{2},...,t_{b}$.  
We take then a dual basis $\Sigma_{j}$ for $H_{2}(E(2,K),\Z)$.

 For sake of comprehension we can assume that the first two elements 
of this dual basis are classes whose geometric representatives are defined in 
the following way (see \cite{FS}, page 387): 
consider the two (fibered) classes $H^{1}(M,\Z)$ of
coordinates $(1,0)$ and $(0,1)$. Their Poincar\'e dual are the classes of
the Seifert surfaces of each link component. 
These have by definition intersection $\delta_{ij}$ with
$\mu(K_{j})$ and are represented by a genus $0$ Riemann surface with 
$1 + 3^{2n}$ discs removed: they span in
fact one link component (the first disk) and intersect the tubular neighborhood of the other in $lk(K_{1},K_{2})$
copies of the meridian (the remaining $3^{2n}$ disks).
We can glue these surfaces to the disc sections of 
$E(1)_{i} \setminus \nu F_{i}$
appearing in the description of $E(2,K)$ (of Equation \ref{ncs})
associated to these two fibered classes.
The closed surfaces $\Sigma_{i}$ arising this way satisfy 
$t_{j}(\Sigma_{i}) = \delta_{ji}$.

 We have now the following simple proposition:

 \begin{pr} \label{prosano}
Let $\varphi = m_{1} \Sigma_{1} + m_{2} \Sigma_{2} + \sum_{j=3}^{b}
m_{j} \Sigma_{j} = (m_{1},m_{2},...,m_{b}) \in H_{2}(E(2,K),\R)$, 
and denote by $(m_{1},m_{2}) \in H^{1}(M,\R)$ the element which has the same two first coordinates: then
\begin{equation} \label{sano} \| \varphi \|_{SW} = \|(m_{1},m_{2})\|_{A} = 
\|(m_{1},m_{2})\|_{T} \end{equation} \end{pr}

\begin{figure}[ht!] 
\centerline{\includegraphics[height=55mm,width=60mm]{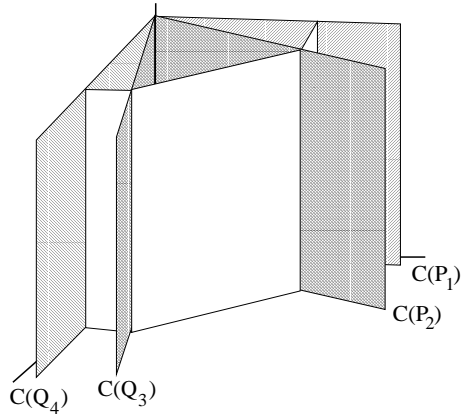}}
\caption{\label{walls}$H_{2}(E(2,K^{(4)}),\R)$ --- the shaded planes
are some of the cylinders over non-fibered rays; the white planes
represents some of the top dimensional faces.}
\end{figure}

\begin{proof} The Seiberg--Witten polynomial of $E(2,K)$ is determined by
the Alex\-ander polynomial of $K$ through Equation \ref{swpol}; therefore, the 
first two coefficients of $\varphi$ determine its SW norm, and that 
coincides with the Alexander norm of $(m_{1},m_{2})$. \end{proof}
We note, as a consequence, that the SW unit ball is a cylinder over the SW 
unit ball on the subspace of $H_{2}(E(2,K),\R)$ spanned by
$\Sigma_{1},\Sigma_{2}$, that we will identify now on with $H^{1}(M,\R)$.
With this identification, the $4n$ top dimensional (dimension $b-1$) faces are
a cylinder over the fibered faces of $H^{1}(M,\R)$; the cone over each top 
dimensional face is bounded by the cylinders over two non-fibered rays, 
cones over lower dimensional (dimension $b-2$) faces.
We will denote these cylinders using the primitive element of the 
non-fibered ray, as
$C(\pm P_{i})$ (or $C(\pm Q_{i}))$); for example, if the
coordinates of $P_{i}$ are $(m_{1},m_{2})$,
\begin{equation} C(P_{i}) = \{ (\mu_{1},\mu_{2},...,\mu_{b}) 
\in H_{2}(E(2,K^{(2n)}),\R) | 
\mu_{1} = \lambda m_{1}, \ \ \mu_{2} = \lambda m_{2}, \ \ \lambda > 0 \}. 
\end{equation}
Figure \ref{walls} illustrates the previous definitions.

 With the knowledge of the SW norm of $E(2,K)$ we are
ready to state our main result. We have: 

\begin{thr} \label{main} {\rm(a)}\qua For each fibered
cohomology classes ${\bf m} \in H^{1}(M,\Z)$ there exist a symplectic 
structure on $E(2,K)$ whose canonical bundle is the image (up to a factor $-1$)
of the Euler class of the fibration under the injective map $H_{1}(M,\Z) 
\rightarrow H^{2}(E(2,K),\Z)$.

 {\rm(b)}\qua Inequivalent fibrations induce inequivalent
symplectic structures on the same smooth manifold $E(2,K)$; in particular, there are $n+1$ inequivalent
symplectic structures on $E(2,K^{(2n)})$. \end{thr}

\begin{proof} (a)\qua For all choice of a fibered class $\bf m$ we
have a symplectic manifold $S^{1} \times M({\bf m})$, with a pair of 
embedded, framed symplectic tori $S^{1} \times
C_{i}$.  $E(1)$, with its standard K\"ahler form, is as well symplectic,
with $F$ an embedded, framed symplectic torus. 
Gompf's theory of symplectic fiber sum (see eg \cite{GS}) guarantees that 
the manifolds defined as in
Equation \ref{ncs} admit a symplectic structure. Proposition \ref{fibsum} tells that 
the smooth manifold carrying the symplectic structure is the same for any 
${\bf m}$. 
Once ${\bf m}$ is fixed, Theorem 3.4 of \cite{McMT} allows to identify the
canonical class of the symplectic structure on $E(2,K)$ as the image, with
reversed sign, of the Euler class of the fibration under the injective map
$ H_{1}(M,\Z) \rightarrow H^{2}(E(2,K),\Z)$ described in Equation \ref{injmap}.
In fact, in \cite{McMT} it is shown that the anticanonical bundle has a 
section whose zero set is homologous to the image of the Euler class 
associated to 
${\bf m}$ cross the class of $S^{1}$. This section is obtained in the 
following way. The anticanonical bundle restricts on $S^{1} \times M$ 
to the pull-back of the plane field defined by the fibration ${\bf m}$, which 
has a preferred section which is inward pointing on the boundary and whose
zero set is the Euler class of ${\bf m}$. This section obviously induces a 
section of the anticanonical restricted to $S^{1} \times M$.
On $E(1) \setminus \nu F$ the 
anticanonical is a trivial bundle (recall that $K_{E(1)} = -F$) and we can 
construct a nowhere vanishing section by restricting a section of $E(1)$ with
a simple zero along $F$: on the boundary $F \times S^{1}$, 
the section can be assumed to depend only on the $S^{1}$ coordinate, 
(as the tangent bundle to the elliptic fiber is canonically trivial), and
outward pointing, because of the simple zero at $F$. 
These two sections glue on the boundary to give a section, 
on $E(2,K)$, whose zero set is the
image of the Euler class under the map of Equation \ref{injmap}.

 Note that the canonical classes arising through this construction satisfy 
(as they should) all Taubes' constraints for $E(2,K)$, 
as they are image 
of a vertex of the polyhedron dual to the Alexander norm ball and 
therefore, by Equation \ref{swpol}, they are the vertices of the convex hull
of basic classes.

To prove part (b) we note that, as symplectic forms with different
canonical bundles are not homotopic, we can construct at least $4n$
non-homotopic symplectic forms on $E(2,K^{(2n)})$ by choosing, in the
symplectic fiber sum definition of Equation \ref{ncs}, elements ${\bf
m}$ laying in different fibered cones of $H^{1}(M^{(2n)},\Z)$ (and
having thus different Euler classes). We need to prove now that the
$4n$ canonical bundles arising this way partition in at least $(n+1)$
orbits under the action of the diffeomorphism group of
$E(2,K^{(2n)})$.

 Diffeomorphisms preserve the convex hull of Seiberg--Witten basic classes, as
well as the unit ball of the Seiberg--Witten norm, which is 
its dual polyhedron. We will study the orbits of the 
vertices of the convex hull (and therefore of the canonical bundles) by
studying the orbits of their dual (top dimensional) face. 
Again, we leave aside the case of $n=1$, where we can easily 
apply divisibility to prove the result.

 Let's assume then that $n \geq 2$. We will study the orbits of the
classes of $H_{2}(E(2,K^{(2n)}),\R)$ laying in the plane spanned by
$\Sigma_{1},\Sigma_{2}$ (which can be identified with
$H^{1}(M^{(2n)},\R)$) with same coordinates as $\pm P_{i}$ or $\pm
Q_{i}$.  For sake of notation, we denote these elements with the same
symbol.  By invariance of the unit ball, diffeomorphisms exchange the
cones $\{ C(\pm P_{j}), C(\pm Q_{n+j}), j=1,...,n\}$: the image of an
element $\pm P_{i}, \pm Q_{i}$ must be therefore a primitive element
laying in one of the $\{C(\pm P_{j}), C(\pm Q_{n+j}),j=1,...,n\}$ and
having the same Seiberg--Witten norm.  Note that there is no reason
why diffeomorphism should preserve the span of
$\Sigma_{1},\Sigma_{2}$. We need to modify a bit the approach to the
proof of Theorem \ref{ineq} to keep track of this.

 Consider one of the four primitive classes $\varphi =
(m_{1},m_{2},0,...,0)$ whose coefficients $m_{1},m_{2}$ lie on
non-fibered rays of $H^{1}(M^{(2n)},\Z)$, and which have the smallest
SW norm. According to Prop \ref{eisneu} and Proposition \ref{prosano},
this must be one of the $\{\pm P_{n}, \pm Q_{n+1}\}$. Take $P_{n}$ for
example.  Under the action of a diffeomorphism $f$, $\varphi$ is sent
to a primitive class $f_{*}\varphi =
(m^{'}_{1},m^{'}_{2},m^{'}_{3},...,m^{'}_{b})$ where, in principle,
$(m_{1}^{'},m_{2}^{'})$ can be divisible, but must be integer multiple
of the coordinates of one of the $\{\pm P_{j},\pm Q_{n+j},j=1,...,n\}$.
Moreover we must have $ \| \varphi \|_{SW} = \| f_{*} \varphi
\|_{SW}$: consequently, by Proposition \ref{prosano} \begin{equation}
\|(m_{1},m_{2})\|_{T} = \|(m_{1}^{'},m_{2}^{'})\|_{T}. \end{equation}
As the norm of $(m_{1},m_{2})$ is minimal, we deduce that in fact
$(m_{1}^{'},m_{2}^{'})$ are the coordinates of one of the $\{\pm
P_{n}, \pm Q_{n+1}\}$. This condition constrains the orbit of the
entire cone $C(P_{n})$; $f_{*} C(P_{n})$ must be one of the four cones
$\{C(\pm P_{n}),$\break$ C(\pm Q_{n+1})\}$. We can repeat this argument for
the other three elements of minimal SW norm, and then pass to the
other four of second minimal norm and so on. Or, more simply, we
observe that as $n \geq 2$ when we know $f_{*}C(P_{n})$ then we can
unambiguously determine the image of the other cones (take $Q_{n+1}$
and proceed as above to determine $f_{*}C(Q_{n+1})$, and then use the
invariance of the unit ball).

 It is easy to check that this result implies that there are at least 
$(n+1)$ orbits of the action of the diffeomorphism group of $E(2,K^{(2n)})$
on the set of $4n$ top dimensional faces of the unit ball, which play the role 
of the fibered faces for the link exterior. 
Again, the possible identifications correspond to reflection of these faces 
about the diagonals $m_{1} = m_{2}$, $m_{1} = - m_{2}$.

 The dual vertices to these faces, and so the canonical bundles, 
lie in at least $(n+1)$ distinguished orbits of the group of diffeomorphism.
This completes the proof of (b). \end{proof}
We finish this Section by pointing out the difference between 
McMullen--Taubes example and ours. In the case of \cite{McMT}, the two
inequivalent symplectic structures arise from inequivalent fibrations, on the
link complement, obtained as restriction of fibrations over the same closed 
three manifold. In particular, the link surgery construction of the manifold
is the same for all fibrations. In our construction, such approach is not 
viable, and we are constrained to use different closed three manifold, and 
different link surgery construction, for each fibration. The extension 
properties of the diffeomorphisms of $E(1)$ allow to obtain, in spite of this,
diffeomorphic manifolds. The nature of our examples, therefore, is in some 
sense more similar to the situation of the ``accidental'' existence of a third 
symplectic structure on the McMullen--Taubes manifold, 
proven in \cite{V} defining a different link surgery 
presentation of the same smooth manifold.

\section{Some examples} \label{n3}
In this Section we will work out some details of the calculation for
$n=2$ and $n=4$, examples that will allow to make some remarks on the 
divisibility of the canonical classes.
We start with $n=2$; the primitive elements of the $4$ non-fibered lines in
$H^{1}(M^{(4)},\Z)$, together with their Thurston (Alexander) norms, are 
\begin{equation} \begin{array}{cc} P_{1} = (3^{3},-1), \ \  \|P_{1}\|_{T} = 
2080 & P_{2} = (3,-1), \ \ \|P_{2} \|_{T} = 256 \\ \\ Q_{3} = (1,-3), \ \ 
\|Q_{3}\|_{T} = 256 & Q_{4} = (1,-3^{3}), \ \ \|Q_{4}\|_{T} = 2080. 
\end{array} \end{equation}
From these data we can deduce, as in Section \ref{classlink}, 
that the three fibered faces contained between $P_{1}$ and $P_{2}$, 
$P_{2}$ and $Q_{3}$, $-Q_{4}$ and $P_{1}$ are not equivalent under 
diffeomorphism (see Figure \ref{piqu}). When we extend this argument to the 
Seiberg--Witten norm on $H_{2}(E(2,K^{(4)}),\R)$, as in Section \ref{defsym}, 
we obtain the proof of the non-equivalence of the top dimensional faces
bounded by $C(P_{1})$ and $C(P_{2})$, $C(P_{2})$ and $C(Q_{3})$, $C(-Q_{4})$ 
and $C(P_{1})$. The canonical classes dual to these faces determine
three inequivalent symplectic structures, as stated in Theorem \ref{main}.

 Let's see how we can obtain the same result with an independent argument:
the knowledge of the non-fibered lines and the Thurston norm allows 
to compute the dual vertex of the fibered face bounded by two non-fibered 
rays. We will determine, more precisely, the vertices of the convex hull of 
the Alexander basic classes. If $(a_{1},a_{2})$ and $(b_{1},b_{2})$ are the
coordinates of two non-fibered primitive elements bounding a fibered cone (eg, $P_{i}$ and $P_{i+1}$, $P_{n}$ and
$Q_{n+1}$ and so on) their Thurston norm is determined by pairing with the 
same dual vertex of coordinates
$(x,y)$, according to the formula \begin{equation} \label{verco}
 \left( \begin{array}{cc}  a_{1} & a_{2} \\ \\ b_{1} & b_{2} 
\end{array} \right)  \left( \begin{array}{c}  x  \\ \\ y  
\end{array} \right) = \frac{1}{2} \left( \begin{array}{c}  \|(a_{1},a_{2})\|_{T} \\ \\ \|(b_{1},b_{2})\|_{T} 
\end{array} \right). \end{equation} From this we can compute the vertices of 
the convex hull of Alexander basic classes: they are
given by the terms \begin{equation} \label{claba} 
t_{1}^{40}t_{2}^{40}, \ \ t_{1}^{38}t_{2}^{-14}, \ \ t_{1}^{32}t_{2}^{-32},
\ \ t_{1}^{14}t_{2}^{-38} \end{equation} and the other $4$ given by Blanchfield duality. 
\begin{figure}[ht!] 
\centerline{\includegraphics[width=60mm, height=40mm]{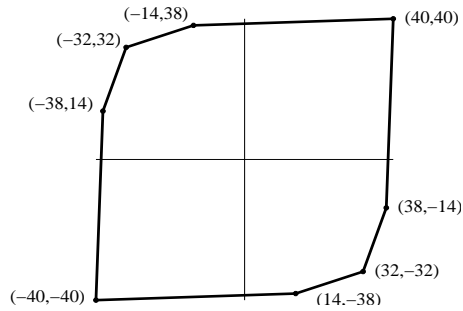}}
\caption{\label{hull}Convex hull of basic classes for $K^{(4)}$} 
\end{figure}

 A check of the
Alexander polynomial (with the help of a computer algebra system) shows that these are in fact the vertices of
the convex hull of the Alexander basic classes, represented in Figure 
\ref{hull}. Note that these vertices are orthogonal to the faces of the 
Thurston unit ball, represented in Figure \ref{piqu}.

 The first three classes of Equation \ref{claba} are dual (up to a factor $2$) 
to the three inequivalent faces. But we can observe, in any case, that they 
have 
different divisibility, and necessarily cannot be joined by a diffeomorphism.
The Euler classes of the fibrations are given by the square 
(in multiplicative notation) of the classes in 
Equation \ref{claba}, and the same is true  (with the aforementioned
identification of the variables
$t_{i}$ as elements of $H^{2}(E(2,K^{(4)}),\Z)$, see Equation \ref{swpol}) 
for the canonical classes corresponding to the inequivalent symplectic structures. In this case,
therefore, the inequivalence of the symplectic structures can be proven by the sole computation of the
Euler classes of the three dimensional fibrations, and in particular without the use of gauge theory.

 The case with $n=4$ appears somewhat different. In that case,
we have $8$ non-fibered lines, and $16$ fibered faces.
The discussion of Section \ref{classlink} leads to the identification of $5$ 
inequivalent classes of fibrations. These classes are represented by the faces 
$(P_{1},P_{2})$, $(P_{2},P_{3})$, $(P_{3},P_{4})$, $(P_{4},Q_{5})$, 
$(-Q_{1},P_{1})$. The vertices of the convex hull of the Alexander basic 
classes dual to these faces can be computed, as before,
using Equation \ref{verco}, and turn out to be
\begin{equation} t_{1}^{3280}t_{2}^{3280},  t_{1}^{3278}t_{2}^{-1094},
 t_{1}^{3272}t_{2}^{-2552},  t_{1}^{3254}t_{2}^{-3038},  
t_{1}^{3200}t_{2}^{-3200}. \end{equation}
Among these classes, the second and the fourth have both divisibility $2$; 
a check of the other divisibilities leads to distinguish at most $4$ classes, 
and consequently at most $4$ symplectic structures on $E(2,K^{(8)})$.
The divisibility criterion for canonical bundles appears therefore weaker 
than the study of the Seiberg--Witten norm, for our family of examples. 
This result should not appear too surprising: we remind the reader 
that all the 
symplectic structures on the McMullen--Taubes manifold defined in \cite{McMT} 
and \cite{V} have the same divisibility.

 It would be interesting to prove that the number of canonical classes of 
different divisibility on $E(2,K^{(2n)})$ grows with no bound with $n$. This 
would give an alternative, non-gauge-theoretic proof of Theorem \ref{prin}.

\end{document}